\documentclass[12pt]{article}                                               
\input{amssym.tex}                                                           
                                                                             
\begin{document}                                                             
%************************************************************                              
\title{Tate curves and UHF-algebras}

\author{Igor  ~Nikolaev
%\footnote{Partially supported 
%by NSERC.}
}

%**************************************************

\date{}
%\date{Nov. 16, 2010}
 \maketitle

%**************************************************

\newtheorem{thm}{Theorem}
\newtheorem{lem}{Lemma}
\newtheorem{dfn}{Definition}
\newtheorem{rmk}{Remark}
\newtheorem{cor}{Corollary}
\newtheorem{prp}{Proposition}
\newtheorem{exm}{Example}
\newtheorem{cnj}{Conjecture}
%*************************************************

\newcommand{\Qcoh}{\hbox{\bf Qcoh}}
\newcommand{\QGr}{\hbox{\bf QGr}}
\newcommand{\ch}{\hbox{\bf Char}}
\newcommand{\jac}{\hbox{\bf Jac}}
\newcommand{\ka}{\hbox{\bf k}}
\newcommand{\n}{\hbox{\bf n}}
\newcommand{\mod}{\hbox{\bf mod}}
\newcommand{\Z}{\hbox{\bf Z}}
\newcommand{\Q}{\hbox{\bf Q}}

%******************************************************************
\begin{abstract}
It is proved that (a stabilization of) the  norm-closure of  a
self-adjoint representation  of the  twisted homogeneous coordinate 
ring of  a  Tate curve  contains  a copy of the  UHF-algebra.

\vspace{7mm}

{\it Key words and phrases: Tate curves,  UHF-algebras}

\vspace{5mm}
{\it MSC: 11G07 (elliptic curves  over local fields);  46L85 (noncommutative topology)}
\end{abstract}

%**************************************************************************
\section{Introduction}
%***************************************************************************
The  {\it Tate curve}   ${\cal E}_q$ is  an affine  cubic over the field  $\Q_p$
of $p$-adic numbers: 
%*************************************************************************
\begin{equation}\label{eq1}
y^2+xy=x^3
-\left(5\sum_{n=1}^{\infty}{n^3q^n\over 1-q^n}\right)x  
-{1\over 12}\sum_{n=1}^{\infty}{(5n^3+7n^5)q^n\over 1-q^n},
\end{equation}
%***************************************************************************
where  $q$ is a $p$-adic integer satisfying condition $0<|q|<1$.   
If    $q^{\Bbb Z}=\{q^n : n\in {\Bbb Z}\}$ is  a lattice
and  $\Q_p^*$ is  the group of units of $\Q_p$,
then   the action $x\mapsto qx$ is discrete;  in particular,   the quotient $\Q_p^*/q^{\Bbb Z}$  
is the Hausdorff topological space.   It was proved by Tate,  that there exists  
an (analytic) isomorphism $\phi: \Q_p^*/q^{\Bbb Z}\to {\cal E}_q$, 
 [Tate 1974]   \cite{Tat1},  p.190.

The  {\it UHF-algebra}  $M_{\ka}$  (``uniformly hyper-finite $C^*$-algebra'')
is  an  inductive limit of the form:
 %*************************************************************************
\begin{equation}\label{eq2}
M_{k_1}({\Bbb C})\to M_{k_1}({\Bbb C})\otimes M_{k_2}({\Bbb C})\to 
M_{k_1}({\Bbb C})\otimes M_{k_2}({\Bbb C})\otimes M_{k_3}({\Bbb C})\to\dots,
\end{equation}
%***************************************************************************
where $M_{k_i}({\Bbb C})$ is a matrix $C^*$-algebra,  $k_i\in\{1, 2, 3,\dots\}$
and $\ka=(k_1,k_2,k_3,\dots)$  [Glimm 1960]  \cite{Gli1}.  Recall that 
the algebras $M_{\ka}$ and $M_{\ka'}$ are said to be   {\it stably isomorphic}
(Morita equivalent)  if  $M_{\ka}\otimes {\cal K}\cong M_{\ka'}\otimes {\cal K}$,
where ${\cal K}$ is the $C^*$-algebra of compact operators;  
the $C^*$-algebra $M_{\ka}\otimes {\cal K}$   is called  a {\it stabilization of }  $M_{\ka}$. 
(The Morita equivalence of $M_{\ka}$ and $M_{\ka'}$ means that  the corresponding 
non-commutative topological spaces are homeomorphic.) 
Glimm classified  stable isomorphism classes of the UHF-algebras 
as follows.    Let $p$ be a prime number and 
$n=\sup~\{0\le j\le \infty : p^j~|~\prod_{i=1}^{\infty}k_i\}$;
denote by $\n=(n_1,n_2,\dots)$ an infinite sequence of $n_i$, where $p_i$
runs through  all primes.  
By ${\Bbb Q}(\n)$ one understands  an additive subgroup of ${\Bbb Q}$
consisting of the rational numbers  whose denominators divide the ``supernatural
number'' $p_1^{n_1}p_2^{n_2}\dots$,   where each $n_j$ belongs to the set
$\{0,1,2,\dots,\infty\}$.   
The ${\Bbb Q}(\n)$ is a dense subgroup
of ${\Bbb Q}$ and every dense subgroup of ${\Bbb Q}$  containing ${\Bbb Z}$ is given by ${\Bbb Q}(\n)$
for some $\n$,   see e.g.  [R\o rdam, Larsen \& Laustsen  2000]  \cite{RLL},  Proposition 7.4.3 (i).  
The UHF-algebra $M_{\ka}$ and the group ${\Bbb Q}(\n)$
are connected by the formula $K_0(M_{\ka})\cong {\Bbb Q}(\n)$, where $K_0(M_{\ka})$
is the $K_0$-group of the $C^*$-algebra  $M_{\ka}$.     
The UHF-algebras  $M_{\ka}$ and $M_{\ka'}$ are stably isomorphic
if and only if  $r{\Bbb Q}(\n)=s{\Bbb Q}(\n')$ for
some positive integers $r$ and $s$,  see  [Effros 1981]   \cite{E}, p. 28 or
[Glimm 1960]  \cite{Gli1}.

The aim of our note is a  link  between  geometry of the Tate curves 
${\cal E}_q$ and  a class of the UHF-algebras  $M_{\ka}$.
Such a link is based on the ideas   of  non-commutative algebraic geometry;
we refer the reader to [Stafford \& van~ den ~Bergh 2001]  \cite{StaVdb1} for 
an excellent introduction to the area. 
Namely,   we consider   the closure of  a  self-adjoint representation of 
the  twisted homogeneous coordinate ring of the Tate curve ${\cal E}_q$ 
by the linear operators on a Hilbert space;  
such a closure  is always a {\it $C^*$-algebra},  see e.g. 
[Murphy 1990]   \cite{M}.  An independent study of 
a related problem has been undertaken  by 
[Carey, Marcolli \& Rennie 2011] \cite{CaMaRe1}. 
Let us pass to an exact statement of our result.

It is known  that each $p$-adic integer $0<|q|<1$ is the limit of a convergent 
series of the rational integers $\alpha_k=\sum_{i=1}^k b_ip^i$,  where 
$0\le b_i\le p-1$ and  $q=\sum_{i=1}^{\infty} b_i p^i$,   see e.g. [Gouv\^ ea 1993]  \cite{G}, p. 66
for the details.
For each $\alpha_k$  one can  define a  supernatural number  $\n(\alpha_k)$
of the form $p_1^{\infty}\dots p_s^{\infty}$,  where 
${\cal P}_k:=\{p_1,\dots, p_s\}$ the finite  set of all primes dividing $\alpha_k$.
By $M_{\alpha_k}$  we shall understand (the stable isomorphism class of)  an  UHF-algebra,
such that $K_0(M_{\alpha_k})\cong {\Bbb Q}(\n(\alpha_k))$. 
Let $\{\pi_1, \pi_2,\dots\}$ be a (finite or infinite)  set of all primes,   such that  
%*************************************************************************
%\begin{equation}\label{eq3}
$\pi_j\in \bigcup_{k=1}^{\infty}  {\cal P}_k$. 
%\end{equation}
%***************************************************************************
By $\n(q)$ we shall understand a supernatural number of the form 
$\pi_1^{\infty} \pi_2^{\infty} \dots$ and by $M_q$ an UHF-algebra,
such that $K_0(M(q))\cong {\Bbb Q}(\n(q))$;  in other words,  $M_q$
is the smallest UHF-algebra containing all the UHF-algebras $M_{\alpha_k}$.   
 Our main result can be stated as follows.    
%************************************************************************
\begin{thm}\label{thm1}
A stabilization of the  norm-closure of  a  self-adjoint  representation of  the (quotient  
 of)  twisted homogeneous coordinate ring of the Tate curve
${\cal E}_q$ contains a copy  the UHF-algebra $M_q$. 
\end{thm}
%*************************************************************************
The note is organized as follows.   Section 2 introduces notation 
and some preliminary facts.   Theorem \ref{thm1} is proved in Section 3;
the proof is based on lemma \ref{lm3}.  
A numerical example illustrating theorem \ref{thm1} is constructed
in Section 4.

%**************************************************************************
\section{Preliminaries}
%***************************************************************************
This section contains a brief description of the twisted  homogeneous  coordinate
ring of an elliptic curve -- the so-called Sklyanin algebra [Sklyanin  1982]  \cite{Skl1};
for a general theory of such rings we refer the reader to  
 [Stafford \& van~ den ~Bergh 2001]  \cite{StaVdb1}. 
The Cuntz-Krieger $C^*$-algebras were introduced 
in [Cuntz \& Krieger  1980 ]   \cite{CuKr1};    their $K$-theory and the crossed product 
structure are discussed   in  [Blackadar  1986]  \cite{B}, Exercise 10.11.9.  
 We consider the two-dimensional  Cuntz-Krieger algebras only;  the latter  link
  the Sklyanin $\ast$-algebras  with  the UHF-algebras.

%**************************************************************************
\subsection{Sklyanin algebras}
%***************************************************************************
Denote  by  $k$ be a field of  characteristic different from  two.  
By a   four-dimensional {\it Sklyanin algebra}   ${\goth S}_{\alpha,\beta,\gamma}(k)$
one understands  a free $k$-algebra  on  four generators  $x_i$  which   satisfy 
six  quadratic relations:
%*******************************************************************
\begin{equation}\label{eq4}
\left\{
\begin{array}{ccc}
x_1x_2-x_2x_1 &=& \alpha(x_3x_4+x_4x_3),\\
x_1x_2+x_2x_1 &=& x_3x_4-x_4x_3,\\
x_1x_3-x_3x_1 &=& \beta(x_4x_2+x_2x_4),\\
x_1x_3+x_3x_1 &=& x_4x_2-x_2x_4,\\
x_1x_4-x_4x_1 &=& \gamma(x_2x_3+x_3x_2),\\ 
x_1x_4+x_4x_1 &=& x_2x_3-x_3x_2,
\end{array}
\right.
\end{equation}
%*****************************************************************  
where $\alpha,\beta,\gamma\in k$ and    $\alpha+\beta+\gamma+\alpha\beta\gamma=0$.
Assume that  $\alpha\not\in \{0;\pm 1\}$;    then  algebra ${\goth S}_{\alpha,\beta,\gamma}(k)$  defines a
non-singular  elliptic  curve   given as  the intersection of two quadrics:
%*************************************************************************
%\begin{equation}
\begin{eqnarray}\label{eq5} 
{\cal E}(k) &= \{(u,v,w,z)\in {\Bbb P}^3(k) ~|~ u^2+v^2+w^2+z^2=\cr
 &\cr
 &=  {1-\alpha\over 1+\beta}v^2+{1+\alpha\over 1-\gamma}w^2+z^2=0\}
\end{eqnarray}
%\end{equation}
%****************************************************************************
together with an automorphism $\sigma: {\cal E}(k)\to {\cal E}(k)$,
see [Sklyanin 1982]  \cite{Skl1} and [Smith \& Stafford 1992]  \cite{SmiSta1}, 
p. 267.   The critical fact  that  we shall use in the future is   the following isomorphism:
%*****************************************************************
\begin{equation}\label{eq6}
\QGr~ ({\goth S}_{\alpha,\beta,\gamma}(k)~/~\Omega)\cong  \Qcoh~({\cal E}(k)),
\end{equation}
%********************************************************************
where $\QGr$ is a category of the quotient graded modules over 
the algebra ${\goth S}_{\alpha,\beta,\gamma}(k)$ modulo torsion, $\Qcoh$ a category 
of the quasi-coherent  sheaves on ${\cal E}(k)$ and $\Omega\subset {\goth S}_{\alpha,\beta,\gamma}(k)$ 
a two-sided ideal generated by the  central  elements
$\Omega_1 = x_1^2+x_2^2+x_3^2+x_4^2$ and 
$\Omega_2 = x_2^2+{1+\beta\over 1-\gamma}x_3^2+{1-\beta\over 1+\alpha}x_4^2$,
see   [Sklyanin 1982]  \cite{Skl1},   Theorem 2. 
Since (\ref{eq6})  coincides  with well-known isomorphism  linking  projective variety 
and  its homogeneous coordinate ring  in the classical  (i.e. commutative) algebraic geometry
but algebra    ${\goth S}_{\alpha,\beta,\gamma}(k)$ is no longer commutative (unless $\sigma$
is trivial),  one calls  ${\goth S}_{\alpha,\beta,\gamma}(k)$ (modulo an ideal) 
a {\it twisted homogeneous coordinate ring}  of elliptic curve ${\cal E}(k)$. 
For a general  theory we refer the reader to  
 [Stafford \& van~ den ~Bergh 2001]  \cite{StaVdb1}.

%***************************************************
\begin{lem}\label{lm-1}
For $\beta=1$ and $\gamma=-1$ elliptic curve given by equations (\ref{eq5})
is isomorphic to such in the Legendre normal form:
%*********************************************************************
\begin{equation}\label{eq7}
y^2=x(x-1)(x-\alpha). 
\end{equation}
%**********************************************************************
 \end{lem}
%***************************************************
{\it Proof.} 
For $\beta=1$  and $\gamma=-1$,  one can write  (\ref{eq5})  in the form:
%*******************************************************************
\begin{equation}\label{eq8}
\left\{
\begin{array}{ccc}
(1-\alpha)v^2+(1+\alpha)w^2+2z^2 &=& 0,\\
u^2+v^2+w^2+z^2&=& 0. 
\end{array}
\right.
\end{equation}
%***************************************************************** 
We shall pass in (\ref{eq8}) from variables $(u,v,w,z)$ to
the new variables $(X,Y,Z,T)$ given by the formulas
%*******************************************************************
\begin{equation}\label{eq9}
\left\{
\begin{array}{ccc}
u^2 &=& T^2,\\
v^2 &=& {1\over 2}Y^2-{1\over 2}Z^2-T^2,\\
w^2 &=& X^2+{1\over 2}Y^2-{1\over 2}Z^2-T^2,\\
z^2 &=& Z^2.
\end{array}
\right.
\end{equation}
%***************************************************************** 
Then equations (\ref{eq8}) take the form  
%*******************************************************************
\begin{equation}\label{eq10}
\left\{
\begin{array}{ccc}
\alpha X^2+Z^2-T^2 &=& 0,\\
X^2+Y^2-T^2 &=& 0. 
\end{array}
\right.
\end{equation}
%***************************************************************** 
Let us consider another (polynomial) transformation
$(x,y)\mapsto (X,Y,Z,T)$ given by the formulas
%*******************************************************************
\begin{equation}\label{eq11}
\left\{
\begin{array}{ccc}
X &=& -2y,\\
Y &=& x^2-1+\alpha,\\
Z &=& x^2+2(1-\alpha) x+1-\alpha,\\
T &=& x^2+2x+1-\alpha.
\end{array}
\right.
\end{equation}
%***************************************************************** 
Then both of  the equations (\ref{eq10}) give us the equation
$y^2=x(x+1)(x+1-\alpha)$,  which after a shift $x'=x+1$ 
takes the Legendre form $y^2=x(x-1)(x-\alpha)$.
Lemma \ref{lm-1} follows.
$\square$

%******************************************************************
\begin{cor}\label{cr1}
Whenever $\beta=1$ and $\gamma=-1$,  
one can replace in (\ref{eq5})  $\alpha$ by $1-\alpha$.
\end{cor}
%*******************************************************************
{\it Proof.}   For elliptic curve ${\cal E}_{\alpha}$ in the Legendre form 
$y^2=x(x-1)(x-\alpha)$ the $j$-invariant is given by well-known formula
$j({\cal E}_{\alpha})=2^6 {(\alpha^2-\alpha+1)^3\over \alpha^2(\alpha-1)^2}$.
It is verified directly,  that  $j({\cal E}_{\alpha})=j({\cal E}_{1-\alpha})$,
i.e.  ${\cal E}_{\alpha}$ is isomorphic to ${\cal E}_{1-\alpha}$.
$\square$

%**************************************************************************
\subsection{Cuntz-Krieger algebras}
%***************************************************************************
Let $A$ be a two-by-two matrix with the non-negative
integer entries $a_{ij}$, such that every row and every column of $A$
is non-zero. 
The two-dimensional {\it Cuntz-Krieger algebra} ${\cal O}_A$
is  a $C^*$-algebra of bounded linear operators on a Hilbert space
${\cal H}$, which is  generated by the partial isometries $s_1$ and $s_2$,
and  relations:
%*******************************************************************
\begin{equation}\label{eq12}
\left\{
\begin{array}{ccc}
s_1^*s_1 &=& a_{11}s_1s_1^*+a_{12}s_2s_2^*,\\
s_2^*s_2 &=& a_{21}s_1s_1^*+a_{22}s_2s_2^*,\\
Id &=& s_1s_1^*+s_2s_2^*,
\end{array}
\right.
\end{equation}
%******************************************************************************  
where $Id$ is the identity operator on ${\cal H}$.  
If one defines $x_1=s_1, x_2=s_1^*,  x_3=s_2$ and  $x_4=s_2^*$, then it is 
easy to see, that   ${\cal O}_A$ contains a dense sub-algebra ${\cal O}_A^0$,
which is a free ${\Bbb C}$-algebra on four generators $x_i$ and three quadratic relations:
%*******************************************************************
\begin{equation}\label{eq13}
\left\{
\begin{array}{ccc}
x_2x_1 &=& a_{11}x_1x_2+a_{12}x_3x_4,\\
x_4x_3 &=& a_{21}x_1x_2+a_{22}x_3x_4,\\
1 &=& x_1x_2+x_3x_4,
\end{array}
\right.
\end{equation}
%******************************************************************************  
and an involution acting by  the formula:
%****************************************************************************
\begin{equation}\label{eq14}
x_1^*=x_2,  \qquad  x_3^*=x_4. 
\end{equation}
%***************************************************************************
Notice,  that equations (\ref{eq13}) are invariant  of  this   involution. 
%******************************************************************************
\begin{lem}\label{lm0}
{\bf (${\cal O}_A$ as a crossed product)}
Let ${\Bbb A}$ be a stationary AF-algebra given by the inductive limit 
${\Bbb Z}^2 \buildrel\rm\ A^T \over\longrightarrow {\Bbb Z}^2
   \buildrel\rm A^T \over\longrightarrow\dots,$
where $A^T$ is the transpose of matrix $A$.  Then 
%*******************************************************************************
\begin{equation}\label{eq15}
{\cal O}_A\otimes {\cal K}\cong {\Bbb A}\rtimes_{\alpha} {\Bbb Z},
\end{equation}
%*******************************************************************************   
where ${\cal K}$ is the $C^*$-algebra of compact operators and $\alpha$
the shift automorphism of ${\Bbb A}$.   In particular,  ${\Bbb A}$
is a  sub-$C^*$-algebra of the stabilization of the  Cuntz-Krieger algebra ${\cal O}_A\otimes {\cal K}$.
\end{lem}
%*******************************************************************************
{\it Proof.}  
We refer the reader to [Effros 1980] \cite{E}, Chapter 6 for the definition of stationary 
AF-algebra, shift automorphism, etc.; see [Blackadar 1986] \cite{B},  Chapter V
for the definition of crossed product $C^*$-algebras.  For a  proof of lemma \ref{lm0},
see [Blackadar 1986] \cite{B}, Exercise 10.11.9.
$\square$

%**************************************************************************
\section{Proof  of theorem \ref{thm1}}
%***************************************************************************
Let the ground field be complex numbers, i.e.  $k={\Bbb C}$.
We shall split the proof in a series of lemmas starting with the following 
elementary
%******************************************************************************
\begin{lem}\label{lm1}
The  ideal of free algebra ${\Bbb C}\langle x_1,x_2,x_3,x_4\rangle$ 
generated by   equations (\ref{eq4})  is  invariant   under   involution (\ref{eq14}), 
if and only if,  $\bar\alpha=\alpha, \beta=1$ and $\gamma=-1$. 
\end{lem}
%*******************************************************************************
{\it Proof.} (i) Let us consider the first two equations (\ref{eq4});
this pair is invariant of involution (\ref{eq14}).  Indeed,  by the rules of composition
for an involution
%*******************************************************************
\begin{equation}\label{eq16}
\left\{
\begin{array}{ccc}
(x_1x_2)^* &= x_2^*x_1^* &= x_1x_2,\\
(x_2x_1)^* &= x_1^*x_2^* &= x_2x_1,\\
(x_3x_4)^* &= x_4^*x_3^* &= x_3x_4,\\
(x_4x_3)^* &= x_3^*x_4^* &= x_4x_3.
\end{array}
\right.
\end{equation}
%*****************************************************************  
Since $\alpha^*=\bar\alpha=\alpha$,  the first two equation (\ref{eq4}) 
remain invariant of  involution (\ref{eq14}).

\smallskip
(ii) Let us consider the middle pair of equations (\ref{eq4});
 by the rules of composition for an involution
%*******************************************************************
\begin{equation}\label{eq17}
\left\{
\begin{array}{ccc}
(x_1x_3)^* &= x_3^*x_1^* &= x_4x_2,\\
(x_3x_1)^* &= x_1^*x_3^* &= x_2x_4,\\
(x_2x_4)^* &= x_4^*x_2^* &= x_3x_1,\\
(x_4x_2)^* &= x_2^*x_4^* &= x_1x_3.
\end{array}
\right.
\end{equation}
%*****************************************************************  
One can apply the involution to the first equation 
$x_1x_3-x_3x_1=\beta(x_4x_2+x_2x_4)$; then one gets
$x_4x_2-x_2x_4=\bar\beta(x_1x_3+x_3x_1)$.  But the second
equation says that $x_1x_3+x_3x_1=x_4x_2-x_2x_4$;  the last  
two equations are compatible if and only if $\bar\beta=1$. 
Thus, $\beta=1$. 

The second equation in involution writes as $x_4x_2+x_2x_4=x_1x_3-x_3x_1$;
the last equation  coincides with  the first equation for $\beta=1$.   

Therefore,  $\beta=1$ is necessary and sufficient for invariance of the middle pair of
equations (\ref{eq4}) with respect to  involution (\ref{eq14}).

\smallskip
(iii) Let us consider the last pair of equations (\ref{eq4});
 by the rules of composition for an involution
%*******************************************************************
\begin{equation}\label{eq18}
\left\{
\begin{array}{ccc}
(x_1x_4)^* &= x_4^*x_1^* &= x_3x_2,\\
(x_4x_1)^* &= x_1^*x_4^* &= x_2x_3,\\
(x_2x_3)^* &= x_3^*x_2^* &= x_4x_1,\\
(x_3x_2)^* &= x_2^*x_3^* &= x_1x_4.
\end{array}
\right.
\end{equation}
%*****************************************************************  
One can apply the involution to the first equation 
$x_1x_4-x_4x_1 = \gamma(x_2x_3+x_3x_2)$; then one gets
$x_3x_2-x_2x_3=\bar\gamma(x_4x_1+x_1x_4)$.  But the second
equation says that $x_1x_4+x_4x_1 = x_2x_3-x_3x_2$;  the last  
two equations are compatible if and only if $\bar\gamma=-1$. 
Thus, $\gamma=-1$. 

The second equation in involution writes as $x_3x_2+x_2x_3=x_4x_1-x_1x_4$;
the last equation  coincides with  the first equation for $\gamma=-1$.   

Therefore,  $\gamma=-1$ is necessary and sufficient for invariance of the last pair of
equations (\ref{eq4}) with respect to  involution (\ref{eq14}).

\smallskip
(iv) It remains to verify that condition $\alpha+\beta+\gamma+\alpha\beta\gamma=0$
is satisfied by $\beta=1$ and $\gamma=-1$ for any $\alpha\in k$.
Lemma \ref{lm1} follows.
$\square$

%***********************************************************************
\begin{rmk}\label{rmk1}
\textnormal{
The Sklyanin algebra ${\goth S}_{\alpha, 1, -1}({\Bbb C})$  with 
$\alpha\in {\Bbb R}$ is a  $\ast$-algebra with the involution 
$x_1^*=x_2$ and $x_3^*=x_4$.
}
\end{rmk}
%*************************************************************************

%***************************************************
\begin{lem}\label{lm2}
The first pair  of equations  (\ref{eq4}):
%*******************************************************************
\begin{equation}\label{eq19}
\left\{
\begin{array}{ccc}
x_1x_2-x_2x_1 &=& \alpha(x_3x_4+x_4x_3),\\
x_1x_2+x_2x_1 &=& x_3x_4-x_4x_3
\end{array}
\right.
\end{equation}
%***************************************************************** 
is equivalent to the pair:
%*******************************************************************
\begin{equation}\label{eq20}
\left\{
\begin{array}{ccc}
x_2x_1 &=& {1+\alpha\over 1-\alpha}x_1x_2-{2\alpha\over 1-\alpha}x_3x_4,\\
x_4x_3 &=& -{2\over 1-\alpha} x_1x_2+ {1+\alpha\over 1-\alpha}x_3x_4.
\end{array}
\right.
\end{equation}
%***************************************************************** 
\end{lem}
%*********************************************************
{\it Proof.} 
Let us isolate $x_2x_1$ and $x_4x_3$ in equations (\ref{eq19});
for that,  we shall write (\ref{eq19}) in the form
%*******************************************************************
\begin{equation}\label{eq21}
\left\{
\begin{array}{ccc}
x_2x_1+\alpha x_4x_3 &=& x_1x_2-\alpha x_3x_4,\\
x_2x_1+x_4x_3 &=& -x_1x_2+x_3x_4.
\end{array}
\right.
\end{equation}
%***************************************************************** 
Consider (\ref{eq21}) as  a linear system of equations relatively
$x_2x_1$ and $x_4x_3$;   since $\alpha\ne 1$,  it  has a unique
solution
%*******************************************************************
\begin{equation}\label{eq22}
\left\{
\begin{array}{ccccc}
x_2x_1 &=& {1\over 1-\alpha} \left|\matrix{x_1x_2-\alpha x_3x_4 & \alpha\cr
                                                                  -x_1x_2+x_3x_4 &  1}\right|
&=&                                                                  
{1+\alpha\over 1-\alpha}x_1x_2-{2\alpha\over 1-\alpha}x_3x_4,\\
&&&&\\
x_4x_3 &=& {1\over 1-\alpha} \left|\matrix{1 & x_1x_2-\alpha x_3x_4 & \cr
                                                                   1 & -x_1x_2+x_3x_4}\right|
&=&
-{2\over 1-\alpha} x_1x_2+{1+\alpha\over 1-\alpha}x_3x_4.
\end{array}
\right.
\end{equation}
%***************************************************************** 
Lemma \ref{lm2} follows.
$\square$

%***********************************************************************
\begin{rmk}\label{rmk2} 
\textnormal{
In new variables $(x_2x_1)'=(1-\alpha)x_2x_1$ and 
$(x_4x_3)'=(1-\alpha)x_4x_3$ the system of equations (\ref{eq20}) 
can be written in the form:
%*******************************************************************
\begin{equation}\label{eq23}
\left\{
\begin{array}{ccc}
x_2x_1 &=& (1+\alpha)x_1x_2- 2\alpha  ~x_3x_4,\\
x_4x_3 &=& - 2  x_1x_2+ (1+\alpha)x_3x_4.
\end{array}
\right.
\end{equation}
%***************************************************************** 
}
\end{rmk}
%*************************************************************************
%***************************************************
\begin{lem}\label{lm3}
{\bf (Main lemma)}
Let  $\alpha_k=\sum_{i=1}^k b_ip^i$ be a rational  integer and $A=\left(\small
\matrix{1+\alpha_k & -2\alpha_k\cr -2 & 1+\alpha_k}\right)$.  
Suppose $M_{\alpha_k}$ is the UHF-algebra  defined in Section 1
and ${\Bbb A}$ an AF-algebra introduced in lemma \ref{lm0}.
Let $I_0$ be the  (two-sided) ideal of the Sklyanin  $\ast$-algebra 
${\goth S}_{\alpha_k, 1, -1}({\Bbb C})$  generated by relation $x_1x_2+x_3x_4=1$
and $J_0$ the ideal of ${\cal O}_A^0$ generated  by  relations
$x_4x_2-x_1x_3=x_3x_1+x_2x_4=x_4x_1-x_2x_3=x_3x_2+x_1x_4=0$.
Then there exists a $\ast$-isomorphism 
%*********************************************************************************
\begin{equation}\label{eq24}
{\goth S}_{\alpha_k, 1, -1}({\Bbb C}) ~/ ~I_0 \cong  {\cal O}_A^0 ~/~ J_0,
\end{equation}
%*******************************************************************************
where 
%*********************************************************************************
\begin{equation}\label{eq25}
\overline{{\cal O}_A^0} \cong {\cal O}_A  \quad  \hbox{and} 
\quad {\cal O}_A\otimes {\cal K} \supset {\Bbb A} \supset M_{\alpha_k}
\end{equation}
%*******************************************************************************
are inclusions of the $C^*$-algebras. 
%***************************************************************** 
\end{lem}
%*********************************************************
{\it Proof.} 
Recall that there exists a dense inclusion ${\Bbb Z}\hookrightarrow \Z_p$ 
given by formula $\alpha_k\mapsto \sum_{i=1}^k b_ip^i$,  where 
$0\le b_i\le p-1$ are integer numbers, see e.g.   [Gouv\^ ea 1993]  \cite{G}, 
Proposition 3.3.4 (ii);  we shall use the inclusion to identify $p$-adic integer
$\sum_{i=1}^k b_ip^i$  with the corresponding rational integer $\alpha_k$.
Because $\alpha_k\in {\Bbb R}$, one gets a Sklyanin $\ast$-algebra 
${\goth S}_{\alpha_k, 1, -1}({\Bbb C})$ with the involution $x_1^*=x_2$
and $x_3^*=x_4$, see remark \ref{rmk1}.   
Notice that the inclusion ${\Bbb Z}\hookrightarrow \Z_p$  induces 
a $\ast$-isomorphism between the following Sklyanin algebras:
%*********************************************************************
\begin{equation}
{\goth S}_{\alpha_k, ~1, ~-1}({\Bbb C})\cong
{\goth S}_{\sum_{i=1}^k b_ip^i, ~1, ~-1}(\Q_p).
\end{equation}
%**********************************************************************
However, the norm convergence in  the above families of Sklyanin
algebras is different;  in what follows, we deal with the Sklyanin 
algebras ${\goth S}_{\alpha_k, 1, -1}({\Bbb C})$ endowed with the
usual  archimedian norm.

\bigskip
{\bf Part  I.}
To prove formula (\ref{eq24}),  one compares relations (\ref{eq4}) 
defining the Sklyanin algebra  ${\goth S}_{\alpha,\beta,\gamma}(k)$ 
with relations (\ref{eq13}) defining
dense sub-algebra ${\cal O}_A^0$ of the Cuntz-Krieger algebra ${\cal O}_A$.
It is easy to see, that ideal $J_0$ is generated by the last four relations of 
system  (\ref{eq4})  corresponding to the case $\beta=-\gamma=1$.    
Likewise,  ideal $I_0$ is generated by the last  relation of system (\ref{eq13}).

As for the first pair of relations of systems (\ref{eq4}) and (\ref{eq13}),
they are identical after a substitution $a_{11}=a_{22}=1+\alpha_k, ~a_{12}=-2\alpha_k$
and $a_{21}=-2$, see also lemma \ref{lm2} and remark \ref{rmk2}.
Thus,  one gets  isomorphism (\ref{eq24}),   where  matrix $A$ is given by
the formula     $A=\left(\small\matrix{1+\alpha_k & -2\alpha_k\cr -2 & 1+\alpha_k}\right)$.

\bigskip
{\bf Part II.}
One can prove inclusions (\ref{eq25}) in the following steps.

\medskip
(i)  The isomorphism  $\overline{{\cal O}_A^0} \cong {\cal O}_A$ follows from
 definition of the Cuntz-Krieger algebra as the norm-closure of algebra
 ${\cal O}_A^0$.
 
 \medskip
 (ii)  In view of lemma \ref{lm0}, there exists a sub-$C^*$-algebra 
 ${\Bbb A}\subset {\cal O}_A\otimes {\cal K}$;   the sub-$C^*$-algebra is the  stationary
 AF-algebra (see [Effros 1980] \cite{E}, Chapter 6) given by the 
 following inductive limit: 
 %*******************************************************************
 \begin{equation}\label{eq26}
 {\Bbb Z}^2 \buildrel\rm
 \left(\small\matrix{1+\alpha_k & -2\cr -2\alpha_k & 1+\alpha_k}\right)
  \over\longrightarrow {\Bbb Z}^2
   \buildrel\rm 
 \left(\small\matrix{1+\alpha_k & -2\cr -2\alpha_k & 1+\alpha_k}\right)  
   \over\longrightarrow\dots
 \end{equation}
 %*********************************************************************
 Notice that since $\alpha_k$ are positive integers, matrix $A^T$ has
 two negative off-diagonal entries.  However,  since $tr~(A^T)>2$
 there exists a matrix in the similarity class of $A^T$ all of whose entries
 are  positive;  the inductive limit (\ref{eq26}) is  invariant of the similarity
 class.  
 
 \medskip
 (iii)    To establish inclusion $M_{\alpha_k}\subset {\Bbb A}$, let us calculate
 the dimension group of  AF-algebra ${\Bbb A}$,   see   [Effros 1980] \cite{E}
 for definition of such a group.  It is known, that  for  stationary AF-algebra
 ${\Bbb A}$ the dimension group is order-isomorphic to ${\Bbb Z}[{1\over\lambda_{A^T}}]$,
 where $\lambda_{A^T}$ is the maximal eigenvalue of matrix $A^T$.   
  We encourage the reader to verify that
  %*********************************************************************************
  \begin{equation}\label{eq27}
  {\Bbb Z} \left[{1\over\lambda_{A^T}}\right]=
  {\Bbb Z}\left[{1+\alpha_k+2\sqrt{\alpha_k}\over (\alpha_k-1)^2}\right]. 
  \end{equation}
  %*******************************************************************************
 It follows from (\ref{eq27}) that
 %*********************************************************************************
  \begin{equation}\label{eq28}
  {\Bbb Z}  \left[{1\over\alpha_k-1}\right] \subset 
  {\Bbb Z} \left[{1\over\lambda_{A^T}}\right]
 \end{equation}
 %*******************************************************************************
 is an inclusion of dimension groups.   In view of corollary \ref{cr1},
 one can replace $\alpha_k-1$ by $-\alpha_k$ in formula (\ref{eq28});
 therefore,   one gets the inclusion:
 %*********************************************************************************
  \begin{equation}\label{eq29}
  {\Bbb Z}\left[{1\over \alpha_k}\right] \subset  
  {\Bbb Z} \left[{1\over\lambda_{A^T}}\right].  
  \end{equation}
  %*******************************************************************************
 Because $\alpha_k\not\in\{0; \pm 1\}$ (see Section 2.1),  one concludes
 that  ${\Bbb Z}\left[{1\over \alpha_k}\right]$ is a dense abelian subgroup
 of the rational numbers ${\Bbb Q}$.  
 It remains to notice that  the dimension group ${\Bbb Z}\left[{1\over \alpha_k}\right]$
is order-isomorphic to such of the UHF-algebra $M_{\alpha_k}$;  
see definition of $M_{\alpha_k}$ in Section 1.    Thus inclusion (\ref{eq29})  
implies the inclusion $M_{\alpha_k}\subset {\Bbb A}$.
Lemma \ref{lm3} follows.
$\square$

\bigskip
Theorem \ref{thm1}  follows from lemma \ref{lm3}. 
$\square$

%**************************************************************************
\section{Example}
%***************************************************************************
We shall consider an example illustrating theorem \ref{thm1}. 
Let $p$ be a prime number  and consider the $p$-adic integer
of the form  $q=p$;  
 notice that  in this case   $b_1=1$ and $b_2=b_3=\dots=0$.
One gets therefore  a supernatural number  $\n(q)$ of the form $p^{\infty}$.
The $\n(q)$   corresponds to a dense subgroup of  ${\Bbb Q}$  
of the  form:
%*************************************************************************
\begin{equation}\label{eq30}
{\Bbb Q}(\n)= {\Bbb Z}\left[{1\over p}\right]. 
\end{equation}
%***************************************************************************
It is easy to see,  that the UHF-algebra corresponding to the Tate curve ${\cal E}_p=\Q_p^*/p^{\Bbb Z}$ 
has the  form:
%*************************************************************************
\begin{equation}\label{eq31}
M_{p^{\infty}}:= M_p({\Bbb C})\otimes M_p({\Bbb C})\otimes\dots
\end{equation}
%***************************************************************************
In virtue  of theorem \ref{thm1},  the UHF-algebra $M_{p^{\infty}}$ is
(a sub-$C^*$-algebra of the stable closure of infinite-dimensional representation of quotient ring of)  
the twisted homogeneous coordinate ring 
of the Tate curve ${\cal E}_p$.    In  particular,  for  ${\cal E}_2$  such a coordinate  ring  is 
the UHF-algebra $M_{2^{\infty}}$;    the latter is known as  
a  Canonical Anticommutation Relations $C^*$-algebra 
(the CAR or Fermion algebra) [Effros 1981]   \cite{E}, p.13.

%**************************************************************************

%**********************************************************

\vskip1cm

\textsc{The Fields Institute for Research in Mathematical Sciences, Toronto, ON, Canada,  
E-mail:} {\sf igor.v.nikolaev@gmail.com}

\end{document}